\documentclass[11pt,a4paper,reqno]{amsart}

\newcommand{\Z}{\mathbb{Z}}
\newcommand{\Q}{\mathbb{Q}}
\newcommand{\N}{\mathbb{N}}

\newtheorem{theorem}{Theorem}
\newtheorem{conjecture}{Conjecture}
\newtheorem{corollary}{Corollary}

\newtheorem{proposition}{Proposition}
\newtheorem{lemma}{Lemma}

\newtheorem{problem}{Problem}

\newcounter{fiddletheoremtemp}
\newenvironment{fiddletheorem}[1]{
    \setcounter{fiddletheoremtemp}{\value{theorem}}
\setcounter{theorem}{#1}\addtocounter{theorem}{-1}
}{
    \setcounter{theorem}{\value{fiddletheoremtemp}}
}

\begin{document}

\title{On Groups and Counter Automata}

\maketitle

\begin{center}

    MURRAY ELDER

    \medskip

    Department of Mathematical Sciences, \ Stevens Institute of Technology, \\
    Hoboken, \ New Jersey 07030, \ USA.

    \medskip

    \texttt{melder@stevens.edu}

    \bigskip

    \bigskip

    MARK KAMBITES

    \medskip

    School of Mathematics, \ University of Manchester, \\
    Manchester M60 1QD, \ England.

    \medskip

    \texttt{Mark.Kambites@manchester.ac.uk} \\

    \bigskip

    \bigskip

    GRETCHEN OSTHEIMER

    \medskip

    Department of Computer Science, \ Hofstra University, \\
    Hempstead, \ New York 11549, \ USA.

    \medskip

    \texttt{cscgzo@hofstra.edu}
\end{center}

\begin{abstract}
We study finitely generated groups whose word problems are accepted by
counter automata.
We show that a group has word problem accepted by a blind $n$-counter
automaton in the sense of Greibach if and only if it is virtually
free abelian of rank $n$; this result, which answers a question of Gilman,
is in a very precise sense an abelian analogue of the Muller-Schupp
theorem. More generally, if $G$ is a virtually abelian group then every
group with word problem recognised by a $G$-automaton is virtually abelian
with growth class bounded above by the growth class of $G$. We consider
also other types of counter automata.
\end{abstract}

\section{Introduction}

Blind counter automata, and the languages they accept, were introduced
and extensively studied by Greibach \cite{Greibach76,Greibach78}. Such an
automaton is a finite state acceptor augmented with a number of integer
counters; these are all initialised to zero, and can be incremented and
decremented during operation, but not read. The automaton accepts a word
exactly if it the reaches an accepting state, with the counters all returned
to zero.

At the same time, an area of considerable interest in combinatorial group
theory is the study of \textit{word problems} of finitely generated groups; it is now
well-known that many structural properties of groups
correspond naturally to language theoretic properties of their word problems,
and vice versa \cite{Boone74,Higman61,Holt05,Muller83}. One of the main objectives of this paper is to give a
complete characterisation of the class of groups whose
word problems are blind counter languages, thus answering a question
posed by Gilman \cite{GilmanWebNotes}.
\begin{theorem}\label{thm_maincounter}
Let $H$ be a finitely generated group. Then the word problem for $H$ is accepted by a blind
$n$-counter automaton if and only if $H$ is virtually free abelian of rank
$n$ or less.
\end{theorem}

As well as being of interest in its own right, Theorem~\ref{thm_maincounter}
forms an important part of a more general framework.
Blind counter automata can be viewed as a class of \textit{$G$-automata}
or \textit{blind group automata}. Recall that a $G$-automaton is a
finite state acceptor augmented with a memory register which stores
an element of a given monoid or group; the automaton cannot read the content of
its register, but it may change the value by multiplying
on the right by some element of the monoid. The register is initialised with
the identity element, and the automaton accepts an input word $w$ exactly if,
by reading this word, it can reach a final state with the register
returned to the identity. Many important classes of languages can be
characterised as $G$-automaton languages for particular groups $G$; 
these include the context-free languages \cite{Chomsky63,Corson05,KambitesGAuto}
and the recursively enumerable languages \cite{Mitrana97}. It follows easily
from the definitions that blind $n$-counter automata are the same thing as
$\Z^n$-automata \cite{KambitesGAuto}; the latter were studied independently
of Greibach's work by Mitrana and Stiebe \cite{Mitrana97}.

Recently, a number of authors have studied connections between the structural
properties of a given group $G$ and of the collection of groups whose word
problems are recognised by $G$-automata \cite{Elder05,Elston04,Gilman98,KambitesGAuto,KambitesGraphRat}.
The second author \cite{KambitesWordProblemsRecognisable}, building upon
previous work of Elston and the third author \cite{Elston04}, has recently given
a complete characterisation of the groups whose word problems are accepted
by \textit{deterministic} $G$-automata, for each group $G$.
\begin{theorem}[Kambites 2006 \cite{KambitesWordProblemsRecognisable}]\label{thm_kambites}
Let $G$ and $H$ be groups, with $H$ finitely generated. Then $H$ has word
problem accepted by a deterministic $G$-automaton if and only if $H$ has a
finite index subgroup which embeds in $G$.
\end{theorem}
An interesting question is that of the extent to which this result can be
extended to the non-deterministic case. A well-known theorem of Muller and
Schupp \cite{Muller83}, when combined with a subsequent result of Dunwoody
\cite{Dunwoody85}, says that a group has context-free word problem if and
only if it is virtually free. If $G$ is a free group of rank $2$ or more
than the $G$-automata accept exactly the context-free languages
\cite{Chomsky63,Corson05,KambitesGAuto}. Moreover, every subgroup of $G$ is
free, and $G$ contains every countable rank free group as a subgroup.
It follows that the result of Muller and Schupp admits the following equivalent
formulation, as a generalisation of Theorem~\ref{thm_kambites} to
non-deterministic automata, in the special case that the register group $G$
is free and non-abelian.
\begin{theorem}[Muller and Schupp 1983 \cite{Muller83}, Dunwoody 1985 \cite{Dunwoody85}]\label{thm_mullerschupp}
Let $G$ be a free group of rank $2$ or more, and $H$ a finitely generated
group. Then the word problem for $H$ is accepted by a $G$-automaton if and
only if $H$ has a finite index subgroup which embeds in $G$.
\end{theorem}

As a consequence of Theorem~\ref{thm_maincounter}, we deduce the following
result, which gives a complete non-deterministic analogue of
Theorem~\ref{thm_kambites} in the case that the register group $G$ is
virtually abelian, and hence also an abelian analogue of
Theorem~\ref{thm_mullerschupp}.
\begin{theorem}\label{thm_maingauto}
Let $G$ be a virtually abelian group and $H$ a finitely generated group.
Then $H$ has word problem accepted by a $G$-automaton if and only if $H$
has a finite index subgroup which embeds in $G$.
\end{theorem}

This generalises and was inspired by the following recent result of Cleary
and the first and third authors \cite{Cleary06}.
\begin{theorem}[Cleary, Elder and Ostheimer 2006 \cite{Cleary06}]\label{thm_cleary}
The word problem for $\Z^p$ is recognised by a $\Z^q$-automaton only
if and only if $p \leq q$.
\end{theorem}
While Theorem~\ref{thm_cleary} was proved directly using only elementary
linear algebra, the more general Theorems~\ref{thm_maincounter} and
\ref{thm_maingauto} seem to require a more complex argument involving
growth of groups and the structure of rational sets in $\Z^n$.

In addition to this introduction, this paper comprises six sections.
Section~\ref{sec_preliminaries} briefly recalls some definitions and
notation which we shall use in the rest of the paper.
Section~\ref{sec_semilinear} obtains some bounding results concerning
minimal elements of intersections of semilinear sets in $\Z^n$; these
may be of some independent interest. In Section~\ref{sec_mainproof},
these bounds are applied to prove that a group with word problem accepted
by a $\Z^n$-automaton has growth bounded above by a polynomial of
degree $n$; it follows by Gromov's polynomial growth theorem
\cite{Gromov81} that such a group is virtually nilpotent. In
Section~\ref{sec_nilpotentabelian} we use a combinatorial result of
Mitrana and Stiebe \cite{Mitrana97} to show that a virtually nilpotent
group with word problem accepted by a $\Z^n$-automaton must in fact be
virtually
abelian. In Section~\ref{sec_conclusion} we combine the results of
Sections~\ref{sec_mainproof} and~\ref{sec_nilpotentabelian} to prove
Theorems~\ref{thm_maincounter} and~\ref{thm_maingauto}. Finally, in
Section~\ref{sec_othercounter}, we consider groups whose word problems
are accepted by other kinds of counter automata.

\section{Preliminaries and Notation}\label{sec_preliminaries}

In this section we recall some basic definitions which we shall need in
the rest of the paper. We assume a basic familiarity with formal languages,
with finite automata over monoids (see \cite{Eilenberg74} for an introduction)
and with elementary group theory (see, for example, \cite{Robinson96}).
Throughout this paper, we write $\N$ for the set of non-negative integers,
including $0$. We apply functions on the right of their arguments, and
compose them from left to right. We denote by the empty word by $\epsilon$.

Let $G$ be a finitely generated group. Recall that a \textit{(finite) choice of
generators} for $G$ is a surjective monoid morphism $\sigma$ from a (finitely
generated) free monoid $X^*$ onto $G$. The choice of generators is called
\textit{symmetric} if $X$ comes equipped with a fixed-point-free involution
$x \mapsto x^{-1}$ such that $(x^{-1}) \sigma = (x \sigma)^{-1}$ for all
$x \in X$; such an involution extends naturally to the whole of $X^*$ by
$(x_1 \dots x_n)^{-1} = x_n^{-1} \dots x_1^{-1}$. We say that a word $w \in X^*$
is a \textit{representative} for the element $w \sigma \in G$. We shall often
leave the morphism $\sigma$ implicit, referring to $X$ as a choice of
generators for $G$. Recall that the \textit{word problem} for $G$
with respect to $\sigma$ is the language of all words $w \in X^*$ such that
$w \sigma = 1$ in $G$.

Now suppose we have fixed a choice of generators $X$ for $G$. For each
$g \in G$, we define the \textit{length} of $g$, denoted $|g|$, to be the
length of the shortest representative for $g$ in $X^*$; if the choice
of generators is
symmetric then we have $|g^{-1}| = |g|$. For $n \in \N$, the \textit{ball
of radius $n$} in $G$ is the set
$$B_n(G) = \lbrace g \in G \mid |g| \leq n \rbrace.$$
Recall that the \textit{growth function} of $G$ (with respect to $X$,
which we usually leave implicit) is
the function
$$\N \to \N, \ n \mapsto | B_n(G) |$$
while the \textit{growth} of $G$ is the asymptotic complexity class of
the growth function. The latter is easily seen to be independent of the
choice of generators.

We now recall the definition of an $M$-automaton.
Let $M$ be a monoid with identity $1$, and $X$ be a finite alphabet. By an
\textit{$M$-automaton $A$ over $X$}, we mean a finite automaton $A$ over
the direct product monoid $M \times X^*$. We say that $A$ \textit{accepts}
a word $w \in X^*$ if it accepts $(1,w) \in M \times X^*$ when considered as a
finite automaton in the usual sense, that is, if there is a path from the
start state to a terminal state labelled $(1,w)$. The \textit{language
recognised} or \textit{accepted} by $A$ is the set of all words $w \in X^*$
which are accepted by $A$. For a more detailed introduction to $M$-automata,
see \cite{Gilman96} or \cite{KambitesGAuto}.

We shall be particularly interested in the case in which the monoid $M$ is
actually a group $G$, and the alphabet $X$ is a (typically symmetric)
monoid generating set for another group $H$. We shall need the following
elementary proposition.
\begin{proposition}\label{prop_subgroup}
Let $G$ be a group, and let $H$ be a finitely generated subgroup of a
finitely generated group
$K$. If the word problem for $K$ (with respect to any choice of generators)
is accepted by a $G$-automaton then the word problem for $H$ (with respect
to any choice of generators) is accepted by a $G$-automaton. If $H$ has
finite index in $K$ then the converse holds.
\end{proposition}
\begin{proof}
Let $\sigma : X^* \to H$ and $\tau : Y^* \to K$ be finite choices of
generators for $H$ and $K$ respectively. For each $x \in X$ choice a
word $w_x \in Y^*$ such that $w_x \tau = x \sigma$, and define a
morphism $\rho : X^* \to Y^*$ by $x \rho = w_x$. It is readily verified
that the word problem for $H$ with respect to $\sigma$ is the inverse
image under $\rho$ of the word problem for $K$ with respect to $\tau$.
Since the class of languages accepted by $G$-automata is closed under
rational transduction \cite[Theorem~6.2]{Gilman96}, and hence under inverse morphism,
this suffices to prove the direct implication.

Conversely, suppose the word problem for $H$ is accepted by a $G$-automaton,
and $H$ has finite index in $K$. By \cite[Theorem~7]{Elston04}, the word
problem for $K$ is accepted by a $H$-automaton, and it follows by
\cite[Corollary~3]{KambitesGAuto} that the word problem for $K$ is accepted
by a $G$-automaton.
\end{proof}
Note that, since a group is a finite index subgroup of itself,
Proposition~\ref{prop_subgroup} implies in particular that recognisability
of the word problem for $H$ by a $G$-automaton is independent of the
choice of generators for $H$.
\begin{corollary}\label{cor_generators}
If the word problem for $H$ with respect to any choice of generators
is accepted by a $G$-automaton, then the word problem for $H$ with respect
to every choice of generators is accepted by a $G$-automaton.
\end{corollary}
We can prove also a similar result for subgroups of the register group.
\begin{proposition}\label{prop_registersubgroup}
Let $G$ be a group, and $H$ a subgroup of $G$. Then any language accepted
by an $H$-automaton is accepted by a $G$-automaton. If $H$ has finite index
in $G$ then the converse holds.
\end{proposition}
\begin{proof}
It follows straight from the definitions that any $H$-automaton can be
viewed as a $G$-automaton accepting the same language, so the direct
implication is immediate.

For the converse, suppose $H$ has finite index in $G$ and that $L$ is
accepted by a $G$-automaton. By \cite[Theorem~7]{Elston04}, the word
problem for $G$ is accepted by an $H$-automaton, and so by
\cite[Corollary~3]{KambitesGAuto}, the language $L$ is accepted by an
$H$-automaton, as required.
\end{proof}

\section{Some Bounds for Linear Maps and Semilinear Sets}\label{sec_semilinear}

The main aim of this section is to establish a bound on the size of the
smallest element of the intersection of two semilinear sets, in terms of
certain parameters of those sets.

Let $n$ be a non-negative integer. We consider the group $\Z^n$ as an
additive subgroup of the row vector space $\Q^n$; hence we use additive
notation. Let $X$ be the set of free generators for $\Z^n$ together with
their inverses, that is, the set of row vectors with a single non-zero component, which takes the
value $1$ or $-1$.
The length function for $\Z^n$ with respect to $X$ is the restriction to
$\Z^n$ of the ``Manhattan taxi'' norm on $\Q^n$, given by
$$\left| (x_1, x_2, \dots, x_n) \right| = |x_1| + |x_2| + \dots + |x_n|$$
where $|x_i|$ for $x_i \in Q$ denotes the absolute value of $x_i$.

Recall that a \textit{linear set} in $\Z^n$ is a set of the form
$$\lbrace s_0 + \lambda_1 s_1 + \dots + \lambda_p s_p \mid \lambda_1, \dots, \lambda_p \in \N \rbrace$$
for some vectors $s_0, s_1, \dots, s_p \in \Z^n$.
A \textit{semilinear set} is a finite union of linear sets. A theorem
of Parikh \cite{Parikh66} states that every context-free subset of $\Z^n$,
and in particular every rational subset of $\Z^n$, is in fact semilinear.

In order to prove the main theorem, we will need some results concerning
linear maps and semilinear sets in $\Z^n$. We begin with the following
elementary linear algebraic proposition 
concerning linear maps.
\begin{proposition}\label{prop_affinebound}
Let $\sigma : \Z^p \to \Z^n$ be a linear map. Then there exist constants
$L, M > 0$
such that if $v \in \Z^n$ lies in $(\N^p) \sigma$, then $v = u \sigma$
for some $u \in \N^p$ with $|u| < L |v| + M$
\end{proposition}
\begin{proof}
Consider first the case in which $\sigma$ is injective. It is readily
verified that $\sigma$ extends uniquely to a linear map $\tau : \Q^p \to \Q^n$ which
is also injective.
Now $\tau$ has an inverse $\tau^{-1} : \mathrm{Im} \ \tau \to \Q^p$, which is
also linear. Hence, there exist
vectors $x_1, \dots, x_n \in \Q^p$ such that
$v \tau^{-1} = v_1 x_1 + \dots + v_n x_n$
for every $v \in \mathrm{Im} \ \tau$, where of course $v_1, \dots, v_n$
denote the components of $v$. Let $L$ be the greatest norm of any of
the vectors $x_i$, and let $M = 0$.

Now suppose $v \in \Z^n$ lies in $\N^p \sigma$. Then letting
$u = v \tau^{-1}$, since $\tau$ is injective, we must have $u \in \N^p$. Now
using the triangle inequality we have:
\begin{align*}
|u|           \ &=|v \tau^{-1}| \\
                &= \ |v_1 x_1 + \dots + v_n x_n| \\
                &\leq \ |v_1| |x_1| + \dots + |v_n| |x_n| \\
                &\leq L (|v_1| + \dots + |v_n|) \\
                &= L |v| + M
\end{align*}
as required.

For the case in which $\sigma$ is not injective, we proceed by induction
on the dimension $p$ of the domain of $\sigma$. First, observe that since
every map with domain $\Z^0$ is injective, the base case $p=0$ is already
proven.

Now let $p \geq 1$, and suppose the claim holds for lesser values of $p$.
It is easily checked that, since $\sigma$ is not
injective, we can choose a non-zero vector $z \in \Z^p$ which lies in
the kernel of $\sigma$; by negating $z$ if necessary, we may assume that at least
one component of $z$ is positive. Let $c$ be the largest positive
component of $z$. Let $x_1, \dots, x_p \in \Z^n$ be vectors such that
$$u \sigma = u_1 x_1 + \dots + u_p x_p$$
for all $u = (u_1, \dots, u_p) \in \Z^p$. For $1 \leq i \leq p$, define
$$\sigma_i : \Z^{p-1} \to \Z^n$$
by currying $\sigma$ with $0$ in the $i$th position, that is, by
$$(x_1, \dots, x_{i-1}, x_{i+1}, \dots, x_p) \sigma_i = (x_1, \dots x_{i-1}, 0, x_{i+1}, \dots, x_p) \sigma.$$
For each $\sigma_i$, let $L_i$ and $M_i$
be the constants given by the inductive hypothesis, and
let $L$ and $M'$ be the greatest of the $L_i$s and $M_i$s respectively.
Let $q$ be the greatest magnitude of any of the vectors $x_i$.

Now suppose $u \sigma = v$ where $u \in \N^p$. Then $(u-z) \sigma = v$, and
either $(u-z) \in \N^p$ or $u_i < c$ for some $i$. By
subtracting $z$ as many times as possible while remaining in $N^p$, we may
assume the latter. Now we have
$$v = u \sigma = u_1 x_1 + \dots + u_p x_p$$
so that
\begin{align*}
v - u_i x_i \ &= \ u_1 x_1 + \dots + u_{i-1} x_{i-1} + u_{i+1} x_{i+1} + \dots + u_p x_p \\
            \ &= \ (u_1, \dots, u_{i-1}, u_{i+1}, \dots, u_p) \sigma_i.
\end{align*}
By the inductive hypothesis, there exists a vector
$$(u_1', \dots, u_{i-1}', u_{i+1}', \dots, u_p') \in \N^{p-1}$$
such that
$$v - u_i x_i \ = \ (u_1', \dots, u_{i-1}', u_{i+1}', \dots, u_p') \sigma_i$$
and
\begin{equation}\label{eqn1}
u_1' + \dots + u_{i-1}' + u_{i+1}' + \dots + u_p' < L_i |v - u_i x_i| + M_i.
\end{equation}
Letting
$$u' \ = \ (u_1', \dots, u_{i-1}', u_i, u_{i+1}', \dots, u_p')$$
we have
$$u' \sigma = (u_1', \dots, u_{i-1}', u_{i+1}', \dots, u_p') \sigma_i + u_i x_i = v.$$
Moreover, by the triangle inequality, we have
$$|v - u_i x_i| \ \leq \ |v| + | - u_i x_i| \ = \ |v| + u_i |x_i|.$$
which combined with \eqref{eqn1} yields
$$u_1' + \dots + u_{i-1}' + u_{i+1}' + \dots + u_p' \ < \ L_i |v| + L_i u_i |x_i| + M_i.$$
Letting $u_i' = u_i < c$ and recalling that $L_i \leq L$, $M_i \leq M'$ and $|x_i| \leq q$ we get
$$|u'| = u_1' + \dots + u_p' \leq L |v| + L c q + M' + c.$$
Since $L$, $M'$, $c$ and $q$ are constants chosen independently of $v$,
setting $M = L c q + M' + c$ suffices to complete the proof.
\end{proof}

We can apply Proposition~\ref{prop_affinebound} to obtain a related bound
on the size of the smallest element of an intersection of linear sets
in $\Z^n$.

\begin{theorem}\label{thm_linearintersect}
Let $n, p, q \geq 0$ and suppose
$s_1, \dots s_p, t_1, \dots t_q \in \Z^n$. Then there exist constants
$P$ and $Q$ such that for any $s_0, t_0 \in \Z^n$ the
linear sets
$$S = \lbrace s_0 + \lambda_1 s_1 + \dots + \lambda_p s_p \mid \lambda_1, \dots, \lambda_p \in \N \rbrace$$
and
$$T = \lbrace t_0 + \mu_1 t_1 + \dots + \mu_q t_q \mid \mu_1, \dots, \mu_q \in \N \rbrace$$
either do not intersect, or have a common element of magnitude less than
$P (|s_0| + |t_0|) + Q$.
\end{theorem}
\begin{proof}
Define a function $\pi : \Z^{p+q} \to \Z^n$ by
$$(\lambda_1, \dots, \lambda_p, \mu_1, \dots, \mu_q) \pi = \lambda_1 s_1 + \dots + \lambda_p s_p - \mu_1 t_1 - \dots - \mu_q t_q.$$
Clearly $\pi$ is linear. Let $L$ and $M$ be the constants given for $\pi$
by Proposition~\ref{prop_affinebound}. Let $c$ be the greatest magnitude
of any of the vectors $s_i$. Let $P = Lc+1$ and $Q = cM$.

Now suppose $S$ and $T$ intersect, and choose $v \in S \cap T$. Then there exist
$\lambda_1, \dots, \lambda_p, \mu_1, \dots, \mu_q \in \N$ such that
$$v \ = \ s_0 + \lambda_1 s_1 + \dots + \lambda_p s_p \ = \ t_0 + \mu_1 t_1 + \dots + \mu_q t_q.$$
Now we have
$$(\lambda_1, \dots, \lambda_p, \mu_1, \dots, \mu_p) \pi = \lambda_1 s_1 + \dots + \lambda_p s_p - \mu_1 t_1 - \dots - \mu_q t_q = t_0 - s_0.$$
so that $t_0 - s_0$ lies in the image of the function $\pi$. Now by Proposition~\ref{prop_affinebound}, there exist 
$\lambda_1', \dots, \lambda_p', \mu_1', \dots, \mu_q' \in \N$ such that
\begin{equation}\label{eqn_lambdabound}
\lambda_1' + \dots + \lambda_p' + \mu_1' + \dots + \mu_q' < L |t_0 - s_0| + M
\end{equation}
and
$$(\lambda_1', \dots, \lambda_p', \mu_1', \dots, \mu_q') \pi = t_0 - s_0.$$
Using the definition of $\pi$, it follows that the element
$$v' = s_0 + \lambda_1' s_1 + \dots + \lambda_p' s_p \ = \ t_0 + \mu_1' t_1 + \dots + \mu_q' t_q.$$
lies in $S \cap T$.
Moreover, using the triangle inequality, we have
\begin{align*}
|v'| \ &= \ |s_0 + \lambda_1' s_1 + \dots + \lambda_p' s_p| \\
       &\leq \ |s_0| + \lambda_1' |s_1| + \dots + \lambda_p' |s_p| \\
       &\leq \ |s_0| + (\lambda_1' + \dots + \lambda_p') c
\end{align*}
Now by \eqref{eqn_lambdabound}, we have
$$\lambda_1' + \dots + \lambda_p' \ \leq \ \lambda_1' + \dots + \lambda_p' + \mu_1' + \dots + \mu_q' \ \leq \ L (|t_0 - s_0|) + M$$
so that
\begin{align*}
|v'| \ &\leq \ |s_0| + Lc |t_0 - s_0| + cM \\
         &\leq \ |s_0| + Lc |t_0| + Lc |s_0| + cM  \ \ \ \text{  (by the triangle inequality)} \\
         &\leq \ (Lc+1) (|s_0| + |t_0|) + cM \\
         &= P (|s_0| + |t_0|) + Q
\end{align*}
as required.
\end{proof}

\begin{corollary}\label{cor_linearintersect}
Let $R$ be a finite set of vectors in $\Z^n$. Then there exist positive
integers $P$ and $Q$ such that any linear sets
$$S = \lbrace s_0 + \lambda_1 s_1 + \dots + \lambda_p s_p \mid \lambda_1, \dots, \lambda_p \in \N \rbrace$$
and
$$T = \lbrace t_0 + \mu_1 t_1 + \dots + \mu_q t_q \mid \mu_1, \dots, \mu_q \in \N \rbrace$$
with $s_1, \dots, s_p, t_1, \dots t_q \in R$ either do not intersect, or have
a common element of magnitude less than
$P (|s_0| + |t_0|) + Q$.
\end{corollary}
\begin{proof}
Let $R'$ denote the (finite) set of all finite sequences of elements from
$R$ which
do not contain the same element twice. For each pair of sequences
$s = (s_1, \dots, s_p) \in R'$ and $t = (t_1, \dots, t_q) \in R'$, let
$P_{st}$ and $Q_{st}$ be the values given by Theorem~\ref{thm_linearintersect},
and choose $P$ and $Q$ to exceed all the $P_{st}$ and $Q_{st}$ respectively.

Now let $S$ and $T$ be as given in the statement, and suppose that they
intersect. Let
$s = (s_1, \dots, s_p)$ and $t = (t_1, \dots, t_q)$. Clearly, we may assume
without loss of generality that $s_i \neq s_j$ for $1 \leq i < j \leq p$
so that $s \in R'$, and similarly $t \in R'$. By Theorem~\ref{thm_linearintersect}, $S$ and $T$ have
a common element of norm less than $P_{st} (|s_0| + |t_0|) + Q_{st}$, which
in turn is less than $P (|s_0| + |t_0|) + Q$, as required.
\end{proof}

We now introduce some terminology for certain parameters of linear and
semilinear sets of $\Z^n$.
Let $p$ be a non-negative integer.
We say that a linear subset $S$ of $\Z^n$ has \textit{constant bound $p$}
if there exist $s_0, \dots, s_q \in G$ with $|s_0| \leq p$ such that 
$$S \ = \ \lbrace s_0 + \lambda_1 s_1 + \dots + \lambda_q s_q \mid \lambda_1, \dots, \lambda_q \in \N \rbrace.$$
Similarly, $S$ has \textit{generator bound $p$} if there exist
$s_0, \dots, s_q \in G$ satisfying the above equation with $|s_i| \leq p$ for $1 \leq i \leq q$.
A semilinear set has constant bound $p$ [respectively, generator bound $p$],
if it is a finite union of linear sets with constant bound $p$ [generator
bound $p$]. We state two elementary properties of these
parameters, which follow directly from the definitions.

\begin{proposition}\label{prop_linearinversebound}
Let $S$ be a linear [semilinear] subset of $\Z^n$ with constant bound
$P$ and generator bound $Q$. Then the set
$$-S = \ \lbrace -s \mid s \in S \rbrace$$
is a linear [semilinear] set with constant bound $P$ and generator
bound $Q$.
\end{proposition}

\begin{proposition}\label{prop_linearsumbound}
Suppose $S$ and $T$ are linear [semilinear] sets with constant bounds $P$
and $Q$ respectively, and generator bounds $R$ and $S$ respectively. Then
the set
$$S+T \ = \ \lbrace s + t \mid s \in S, t \in T \rbrace$$
is a linear [semilinear] set with constant bound $P+Q$ and generator
bound $\max(R,S)$.
\end{proposition}

\begin{theorem}\label{thm_simplerintersect}
Let $K$ be a positive integer. Then there exist positive integers $C$ and
$D$ such that for any semilinear sets $S$ and $T$ in $\Z^n$ with constant
bound $m$ and generator bound $K$, either $S$ and $T$ do not intersect, or
they have a common element of magnitude less than $C m + D$.
\end{theorem}
\begin{proof}
Let $R$ be the (finite) set of all vectors in $\Z^n$ of norm less than $K$,
and let $P$ and $Q$ be the values given by Corollary~\ref{cor_linearintersect}.
Let $C = 2P$ and $D = Q$.

Now suppose $S$ and $T$ are intersecting semilinear sets with constant bound
$m$ and generator bound $K$. Then each is a finite union of linear components
with constant bound $m$ and generator bound $K$, where some component $S'$ of $S$
intersects with some component $T'$ of $T$. Writing
$$S' = \lbrace s_0 + \lambda_1 s_1 + \dots + \lambda_p s_p \mid \lambda_1, \dots, \lambda_p \in \N \rbrace$$
and
$$T' = \lbrace t_0 + \mu_1 t_1 + \dots + \mu_q t_q \mid \mu_1, \dots, \mu_q \in \N \rbrace,$$
by Corollary~\ref{cor_linearintersect} there is an element in $S' \cap T'$,
and hence in $S \cap T$, with norm less than
$$P(|s_0| + |t_0|) + Q \ \leq \ P (m + m) + Q \ = \ 2Pm + Q \ = \ Cm+D.$$
\end{proof}

\section{The Growth of Counter Groups}\label{sec_mainproof}

In this section, we apply the main result of Section~\ref{sec_semilinear}
to show that a group with word problem accepted by a $\Z^n$-automaton
has growth bounded above by a polynomial of degree $n$. It follows by
Gromov's polynomial growth theorem \cite{Gromov81} that such a group is
virtually nilpotent; in Sections~\ref{sec_nilpotentabelian} and \ref{sec_conclusion}
we shall see that every such group is in fact virtually abelian.

\begin{theorem}\label{thm_countergrowth}
Let $H$ be a group. If the word problem for $H$ is accepted by a
$\Z^n$-automaton then $H$ has growth bounded above by a polynomial of
degree $n$.
\end{theorem}
\begin{proof}
Fix a finite, symmetric choice of generators $X$ for $H$; then by
Corollary~\ref{cor_generators} we may let $A$ be a $\Z^n$-automaton
accepting the word problem
$W_X(H)$ of $H$ with respect to $X$. Clearly, we may assume without loss
of generality that every edge in $A$ has label of the form $(g, x)$ for some
$g \in \Z^n$ and $x \in X \cup \lbrace \epsilon \rbrace$.

To prove the lemma, it will suffice to demonstrate the existence of
a function $\sigma : H \to \Z^n$ (which need not be a morphism) and
constants $P$, $Q$ and $R$ such that:
\begin{itemize}
\item[(i)] the pre-image $g \sigma^{-1}$ of each element $g \in \Z^n$ contains
at most $R$ elements; and
\item[(ii)] for any $h \in H$, we have $|h \sigma| \leq P |h| + Q$.
\end{itemize}
Indeed, suppose $\sigma$, $P$, $Q$ and $R$ satisfy these properties.
Then clearly,
$|B_n(H)| \leq R \left( |B_{Pn+Q}(\Z^n)| \right)$ for all $n$. Since the
growth function of $\Z^n$ is a polynomial of degree $n$, and the composition
of a polynomial of degree $n$ with a linear function is a polynomial of
degree $n$, this suffices to show that the growth of $H$ is bounded above
by a polynomial of degree $n$, as required. The rest of this section, then,
is concerned with finding a function $\sigma$ and constants $P$, $Q$ and
$R$ with the above properties.

We begin with a definition. For each pair of states $p, q$ in the automaton $A$, let
$R_{pq}$ denote the set of all element $g \in \Z^n$ such that
$(g, \epsilon)$ labels a path from $p$ to $q$.
By removing from $A$ all edges with label of the form $(g, w)$ for
$w \neq \epsilon$ and then ignoring the second components of the edge
labels, we can clearly obtain a finite automaton over $\Z^n$ which,
with initial state $p$ and terminal state $q$ accepts $R_{pq}$. Thus,
each $R_{pq}$ is a rational subset of $\Z^n$, and so by Parikh's theorem
\cite{Parikh66} is semilinear.

Since there are finitely many sets of the form $R_{pq}$, there clearly exist constants
$F$ and $K$ such that every $R_{pq}$ has constant bound $F$ and generator
bound $K$. Moreover,
since $A$ has only finitely many edges, we can choose $F$ also to be
greater than any $|g|$ such that $A$ has an edge with label of the form
$(g, x)$. Let $C$ and $D$ be the constants given by
Theorem~\ref{thm_simplerintersect} for our chosen value of $K$.
Let $P = 2CF$ and $Q = CF + D$.
Let $R$ be the number of states in the automaton $A$.

The main work of the proof is done in three lemmas.

\begin{lemma}\label{lemma_linearset}
Let $w \in X^*$, and suppose the automaton $A$ has a path from $p$ to $q$
labelled
$(g, w)$. Then there exists a subset $S \subseteq \Z^n$ such that
\begin{itemize}
\item[(i)] $S$ is linear with constant bound $(2 |w| + 1) F$ and generator
bound $K$;
\item[(ii)] $g \in S$;
\item[(iii)] for every $s \in S$, $(s, w)$ labels a path from $p$ to $q$ in $A$.
\end{itemize}
\end{lemma}
\begin{proof}
Suppose $w = w_1 \dots w_m$ with each $w_i \in X$, so that $m = |w|$.
Let $\pi$ be a path from $p$ to $q$ labelled $(g,w)$.
It follows from our assumption above about the edge labels in $A$ that 
the path $\pi$ can be decomposed as
$$\pi \ = \ p_0 e_1 p_1 e_2 \dots e_m p_m$$
where
\begin{itemize}
\item each $p_i$ is a path with label $(k_i, \epsilon)$ for some $k_i \in \Z^n$, and
\item each $e_i$ is an edge with label $(g_i, w_i)$ for some $g_i \in \Z^n$.
\end{itemize}
where of course $k_0 + g_1 + k_1 + \dots + g_m + k_m = g$.

We define
$$T = \lbrace g_1 \rbrace + \dots + \lbrace g_m \rbrace + R_{p_0 \alpha, p_0 \omega} + \dots + R_{p_m \alpha, p_m \omega}$$
where $p_i \alpha$ and $p_i \omega$ denote respectively the start state and
the end state of the path $p_i$.
It is immediate from the definitions that each element $k_i$ lies
in the semilinear set $R_{p_i \alpha, p_i \omega}$. Hence, we see that
$g = k_0 + g_1 + k_1 + \dots + g_n + k_n$ lies in $T$.

Notice that $T$ is the sum of $2m+1$ semilinear sets each of which has
constant bound $F$ and generator bound $K$. It follows by
Proposition~\ref{prop_linearsumbound} that 
$T$ is a semilinear set with constant bound $(2m+1)F$ and generator bound
$K$.
But this means that $T$ is a finite union of linear sets with constant
bound $(2m+1)F$ and generator bound $K$, at least one of which must contain
$k$; choose
$S$ to be such a set. Finally, it follows immediately from the definition
of the sets $R_{pq}$ that for every element of $s \in T$, and hence certainly
for every element $s \in S$, we have that $(s, w)$ labels a path from $p$
to $q$, as required.
\end{proof}

\begin{lemma}\label{lemma_shortword}
Let $w \in X^*$. Then there exists a state $q$ in $A$ and an element
$g \in G$ with $|g| < P |w| + Q$ such that $(g,w)$ labels a
path from the initial state to $q$, and $(-g, w^{-1})$ labels a path
from $q$ to some terminal state.
\end{lemma}
\begin{proof}
Certainly $w w^{-1} \in X^*$ lies in the word problem of $H$ and so is
accepted by the $\Z^n$-automaton $A$. It follows that there is an element
$k \in \Z^n$ and a vertex $q$ in $A$, such that $A$ has
\begin{itemize}
\item a path $\pi_1$ from the initial state to $q$ labelled $(k,w)$; and
\item a path $\pi_2$ from $q$ to a terminal state labelled $(-k,w^{-1})$.
\end{itemize}
By Lemma~\ref{lemma_linearset}, there exist linear sets $S$ and $T$, each
with constant bound $(2|w| + 1) F = (2|w^{-1}| + 1) F$ and generator bound $K$,
such that
\begin{itemize}
\item $k \in S$;
\item $-k \in T$;
\item for every $s \in S$, there is a path from the initial state to $q$ labelled $(s,w)$;
\item for every $t \in T$, there is a path from $q$ to a terminal state labelled $(t,w^{-1})$;
\end{itemize}
By Proposition~\ref{prop_linearinversebound}, the set $-T$ is also linear
with constant bound $(2|w| + 1) F$ and generator bound $K$. Moreover,
$k \in S \cap (-T)$. Hence, by Theorem~\ref{thm_simplerintersect} and
the choice of $C$ and $D$ above, there exists an element $g \in S \cap (-T)$
with
$$|g| \ < \ C(2|w|+1)F+D \ = \ 2CF|w| + CF + D \ = \ P |w| + Q.$$
Now $g \in S$ and $-g \in T$, so that $g$ has the
required properties.
\end{proof}

\begin{lemma}\label{lemma_pigeonhole}
Let $g \in \Z^n$, and let $S \subseteq X^*$ be such that for every $w \in S$,
there exists a state $q_w$ in $A$ such that $(g,w)$ labels a path from the
initial state to $q_w$, and $(-g, w^{-1})$ labels a path from $q_w$ to a
terminal state.
Then $S$ contains representatives for at most $R$ elements of $H$.
\end{lemma}
\begin{proof}
Suppose for a contradiction that $S$ has the given properties, but contains
representatives for
strictly more than $R$ elements of $H$. Then by the pigeon hole principle,
since $A$ has only $R$ states, there exist words $u,v \in S$ such that
$u$ and $v$ represent distinct elements of $H$, but $q_u = q_v$. But
now there is a path from the initial state to $q_u$ labelled $(g,u)$, and
a path from $q_u = q_v$ to a terminal state labelled $(-g, v^{-1})$.
Hence, the $\Z^n$-automaton $A$ accepts the word $u v^{-1}$, so we must have
that $u v^{-1}$ lies in the word problem for $H$. But this means that $u$
and $v$ represent the same element of $H$, giving the desired contradiction.
\end{proof}

We are now ready to define our function $\sigma : H \to \Z^n$.
For each $h \in H$, we choose some word $w \in X^*$ of minimal length
representing $h$. Lemma~\ref{lemma_shortword} guarantees the existence of
at least one element $g$ with $|g| < P |w| + Q = P |h| + Q$ such that
$(g,w)$ labels a path from an initial state
to some state $q$, and $(-g, w^{-1})$ labels a path from $q$ to a terminal
state. We choose $h \sigma$ to be one such element. It is immediate that $\sigma$
satisfies the desired bound restriction.
Moreover, by Lemma~\ref{lemma_pigeonhole}, no more than $R$ choices of $h$
can give rise to the same $h \sigma$; indeed otherwise setting $S$ to be the
set of all the words $w$ chosen for such values of $h$ would contradict
Lemma~\ref{lemma_pigeonhole}. This completes the proof of Theorem~\ref{thm_countergrowth}.
\end{proof}

\section{From Virtually Nilpotent to Virtually Abelian}\label{sec_nilpotentabelian}

The objective of this section is to show that a group which is virtually
nilpotent and which has word problem recognised by a blind counter automaton
is in fact virtually abelian. Our main tool is the following ``interchange
lemma'' of Mitrana and Stiebe \cite{Mitrana97}.

\begin{lemma}[Mitrana and Stiebe 1997 \cite{Mitrana97}]\label{lemma_interchange}
Let $G$ be an abelian group and $L$ a language. Suppose $L$ is accepted by
a $G$-automaton. Then there exists a constant $p$ such that for any word
$w\in L$ of
length at least $p$, and for any factorisation
$$w = v_1 w_1 v_2 w_2 \ldots w_p v_{p+1}$$
with each $|w_i|\geq 1$, there exists integers $r$ and $s$ with
$1 \leq r < s \leq p$ such that the word obtained from
$w$ by interchanging the factors $w_r$ and $w_s$ lies in $L$.
\end{lemma}

We now apply Lemma~\ref{lemma_interchange} to prove the main result of this
section.
\begin{theorem}\label{thm_nilpotentabelian}
Suppose $H$ is a finitely generated virtually nilpotent group with word
problem accepted by a $\Z^n$-automaton. Then $H$ is virtually abelian.
\end{theorem}
\begin{proof}
We claim first that it suffices to prove the result for the case that
$H$ is nilpotent and torsion-free. Indeed, if $H$ is virtually nilpotent
then let $K$ be a nilpotent subgroup of finite index in $H$. Now
by \cite[Theorem~5.5]{Segal83} $K$ embeds in $GL(n, \Z)$ for some
$n$; now by \cite[Lemma~9]{Dixon71} $K$ has a finite index subgroup
$L$ which is torsion-free. Since a subgroup of a nilpotent group is
nilpotent, $L$ is also nilpotent. By Proposition~\ref{prop_subgroup},
the word problem for $L$ is also accepted by a $\Z^n$-automaton. If
the result holds for torsion-free nilpotent groups we may deduce that
$L$ is virtually abelian and hence, since $L$ has finite index in $H$,
that $H$ is virtually abelian.

Suppose, then, that $H$ is not abelian, but is torsion free nilpotent
with centre $Z$. Choose an element $a$ in the second term of the upper
central series for $H$, that is, such that $a$ is not central in $H$ but
$aZ$ is central in $H / Z$. Let $b \in H$ be an element which does not
commute with $a$. Let $c = [a,b] = a^{-1} b^{-1} a b$. Notice that $c$ is
central; indeed $ab Z = aZ bZ = bZ aZ = ba Z$ so that $ab = bac \in baZ$,
from which it follows that $c \in Z$.

Choose a monoid generating set $X$ for $H$ with elements $x$, $y$ and $z$
representing to $a$, $b$ and $c$ respectively. By Corollary~\ref{cor_generators}, the word problem for $H$ with respect
to $X$ is accepted by a $\Z^n$-automaton. Let $p$ be the interchange constant for $A_G$
posited by Lemma~\ref{lemma_interchange}.

Let $w = w_1 \dots w_m \in \lbrace x, y \rbrace^*$.
We define the \textit{exchange index} of $w$ to be the number of distinct
pairs $(i,j)$ with $1 \leq i < j \leq m$ such that $w_i = y$ and $w_j = x$;
it is the minimum number of times one would have to replace a factor $yx$
by $xy$ in order to obtain a word of the form $x^k y^l$.
\begin{lemma}\label{lemma_exchange}
Let $u, v \in \lbrace x, y \rbrace^*$ be words each containing $i$ copies
of the letter $x$, and $j$ copies of the letter $y$. If $u$ and $v$ represent
the same element of $H$, then they have the same exchange index.
\end{lemma}
\begin{proof}
Suppose $u$ and $v$ have exchange indices $p$ and $q$ respectively.
Since the commutator $c$ of $a$ and $b$ is central, we have
$$gabh = gbach = (gbah)c$$
for all $g, h \in H$. It follows by an easy inductive argument that 
$u$ represents the element $a^i b^j c^p$ while $v$ represents the
element $a^i b^j c^q$. If $u = v$ then we have $c^p = c^q$ and hence,
since $H$ is torsion free, $p = q$ as required.
\end{proof}

Consider now the word
$$t_1 = xyxy^{2}xy^{3}\ldots xy^{n},$$ and let $w$ be any word
for its inverse, so that
$t_1 w$ lies in the word problem of $H$.
Choose $n>p$ and set $v_i =x , w_i = y^{i}$
 for $1\leq i\leq n$ and $v_{n+1}=u$.
 Then by Lemma~\ref{lemma_interchange} for some $r$ and $s$ we have
that
$$xyxy^{2}xy^{3}\ldots xy^{s}\ldots xy^{r}\ldots xy^{n}w$$
is also in the word problem for $H$, from which it follows that
$t_1$ and
$$t_2 = xyxy^{2}xy^{3}\ldots xy^{s}\ldots xy^{r}\ldots xy^{n}$$
both represent the same element of $H$. But $t_1$ and $t_2$ contain
the same number of $x$s and $y$s respectively, while a simple
calculation shows that their exchange indices differ by $(r-s)^2$. This
contradicts Lemma~\ref{lemma_exchange}, and hence completes the proof.
\end{proof}

\section{Proofs of the Main Theorems}\label{sec_conclusion}

Combining the results of Sections~\ref{sec_mainproof} and
~\ref{sec_nilpotentabelian} yields our first main result.
\begin{fiddletheorem}{1}
\begin{theorem}
Let $H$ be a finitely generated group. Then the word problem for $H$ is
accepted by a blind $n$-counter automaton if and only if $H$ is virtually
free abelian of rank $n$ or less.
\end{theorem}
\end{fiddletheorem}
\begin{proof}
Suppose the word problem for $H$ is accepted by a blind $n$-counter automaton,
that is, a $\Z^n$-automaton. By Theorem~\ref{thm_countergrowth}, $H$ has
growth bounded above by a polynomial of degree $n$. In particular, $H$ has
polynomial growth, and so by a famous theorem of Gromov \cite{Gromov81}
is virtually nilpotent. It follows by Theorem~\ref{thm_nilpotentabelian}
that $H$ is virtually abelian. It follows easily from the structure theorem
for finitely generated abelian groups \cite[Theorem~4.2.10]{Robinson96}
that $H$ has a finite index subgroup isomorphic to $\Z^k$ for some $k \in \N$.
Moreover, the asymptotic growth of $H$ is a polynomial of degree $k$, so we
must have $k \leq n$.

Conversely, if $H$ is virtually free abelian of rank $n$ then it has a
finite index subgroup isomorphic to $\Z^n$; now by \cite[Theorem~7]{Elston04},
the word problem for $H$ is accepted by a $\Z^n$-automaton, that is, a blind
$n$-counter automaton, as required.
\end{proof}

Theorem~\ref{thm_maincounter} leads to our second main result.
\begin{fiddletheorem}{4}
\begin{theorem}
Let $G$ be a virtually abelian group. Then a group $H$ has word problem
accepted by a $G$-automaton if and only if $H$ has a finite index
subgroup which embeds in $G$.
\end{theorem}
\end{fiddletheorem}
\begin{proof}
Suppose first that the word problem for $H$ is accepted by a
$G$-automaton. Since $G$ has a finite index abelian subgroup, it
follows easily from the structure theorem for finitely generated 
abelian groups \cite[Theorem~4.2.10]{Robinson96} that $G$ has a
finite index subgroup isomorphic to $\Z^n$ for some $n \in \N$. Now
by Proposition~\ref{prop_registersubgroup}, the word problem for $H$ is accepted by a
$\Z^n$-automaton, so by Theorem~\ref{thm_maincounter}, 
$H$ has a finite index subgroup isomorphic to $\Z^k$ for some $k \leq n$.
Now $\Z^k$ embeds in $\Z^n$, which embeds in $G$ as required.

Conversely, if $H$ has a finite index subgroup which embeds in $G$
then by \cite[Theorem~7]{Elston04}, the word problem for $H$ is
accepted by a $G$-automaton.
\end{proof}

We also obtain an answer to a question of the second author
\cite[Question~10]{KambitesWordProblemsRecognisable} for the
class of virtually abelian groups.
\begin{corollary}\label{cor_nondetimpdet}
Let $G$ be a virtually abelian group. Then a group word problem is
accepted by a deterministic $G$-automaton if and only if it is
accepted by a non-deterministic $G$-automaton.
\end{corollary}
\begin{proof}
Suppose $H$ has word problem accepted by a non-deterministic $G$-automaton.
Then by Theorem~\ref{thm_maingauto}, $H$ has a finite index subgroup
which embeds in $G$. Now by \cite[Theorem~7]{Elston04}, the word problem
for $H$ is recognised by a deterministic $G$-automaton. The
converse implication is immediate.
\end{proof}

\section{Other Counter Automata}\label{sec_othercounter}

In this section, we discuss the relationships between blind counter
automata, that is, $\Z^n$-automata, and other notions of counter
automata.

First, we consider the difference between blind and \textit{non-blind}
counter automata. The latter differ from the automata we have studied in
that they are endowed with the extra ability to test whether the value of
a register is zero; for a formal definition see
\cite{Herbst91}, \cite{Herbst93} or \cite[Chapter~7]{Hopcroft69}. It is
well known that a non-blind automaton with two or more counters can simulate
an arbitrary Turing machine \cite[Theorem~7.9]{Hopcroft69}. Since, for a
counter automaton to accept a language, we do not require that the machine should ``terminate''
on words not contained in the language, it follows that $2$-counter non-blind
automata can accept all the recursively enumerable languages and hence, 
by Higman's embedding theorem \cite{Higman61}, word problems of exactly
the finitely presented groups and their finitely generated subgroups.

The remaining case of $1$-counter non-blind automata was studied by
Herbst \cite{Herbst91}; he showed that a group has
word problem accepted by such an automaton if and only if it is virtually
cyclic (see also \cite{Herbst93}). Combining with Theorem~\ref{thm_maingauto},
we see that, when attention is restricted to word problems for groups, there
is no difference in accepting power between blind and non-blind $1$-counter
automata.
\begin{corollary}\label{cor_nonblindimpblind}
A group word problem is accepted by a non-blind $1$-counter automaton
if and only if it is accepted by a blind $1$-counter automaton, that
is, a $\Z$-automaton.
\end{corollary}
The proof of Herbst's result \cite{Herbst91}, and hence that of
Corollary~\ref{cor_nonblindimpblind}, depends upon the
Muller-Schupp theorem \cite{Muller83} and thus also
on Stallings' theory of ends \cite{Stallings68} and the accessibility of
finitely presented groups \cite{Dunwoody85}. One is drawn is ask if there
is an easier approach.
\begin{problem}
Find an elementary proof of Corollary~\ref{cor_nonblindimpblind}.
\end{problem}
Such a proof combined with our proof of Theorem~\ref{thm_maincounter}
above would yield also a purely combinatorial proof of the fact that every
finitely generated
group with word problem accepted by a non-blind $1$-counter automaton has
linear growth. Since there is also an elementary proof that the only finitely
generated groups of linear growth are the virtually cyclic groups
\cite{Wilkie84}, it would also give a straightforward combinatorial proof
of the result of Herbst itself.

Another notion of counter automaton has been studied by Cho \cite{Cho06}. We
use here a definition easily seen to be equivalent in accepting power to his,
although our notation is chosen to make clearer the relationship between
his automata and $G$-automata. A
\textit{$k$-counter Cho automaton} over an alphabet $X$ consists of
\begin{itemize}
\item a finite automaton over the direct product monoid $\N^k \times X^*$,
with edge labels drawn from $\N^k \times X$; and
\item for each state $v$ in the automaton, a semilinear subset $S_v$ of
$\N^k$.
\end{itemize}
A word $w \in X^*$ is \textit{accepted} by the automaton if there
is a path from the initial state to some state $v$ with label $(g,w)$, such
that $g \in S_v$. The automaton is called \textit{deterministic} if for
each state $p$ and letter $x \in X$, there is at most one edge leaving $p$
with label of the form $(g, x)$.
For a detailed discussion of these automata, and the languages they accept,
see \cite{Cho06}. The following theorem characterises the groups whose word
problems are accepted by Cho automata.
\begin{theorem}\label{thm_cho}
Let $H$ be a finitely generated group. Then the following are equivalent
\begin{itemize}
\item[(i)] the word problem for $H$ is accepted by a deterministic Cho automaton;
\item[(ii)] the word problem for $H$ is accepted by a Cho automaton;
\item[(iii)] $H$ is virtually abelian.
\end{itemize}
\end{theorem}
\begin{proof}
That (i) implies (ii) is immediate, so suppose that (ii) holds, and
let $A$ be a $k$-counter Cho automaton accepting the word problem for
$H$. We view $\N^k$ as embedded in the natural way into $\Z^k$, so that
$A$ can be viewed as a $\Z^k$-automaton. We construct from $A$ a new
$\Z^k$-automaton by adding some extra states and edges as follows.

For each state $p$, write the semilinear set $S_p$ as a union
$$S_p \ = \ S_p^1 \cup S_p^2 \cup \dots \cup S_p^j$$
where each set $S_p^i$ is linear. For each $i$, let
$v_0^i, v_1^i, \dots, v_m^i$ be such that
$$S_p^i \ = \ \lbrace v_0^i + \lambda_1 v_1^i + \dots + \lambda_m v_m^i \mid \lambda_1, \dots, \lambda_m \in \N \rbrace$$
For each $i$, the new automaton $B$ has an extra state $p_i$, an edge from
$p$ to $p_i$ labelled $(-v_0^i, \epsilon)$, and edges from $p_i$ to $p_i$
labelled $(-v_q^i, \epsilon)$ for $1 \leq q \leq m$.
The initial state of $B$ is the initial state of $A$, while the terminal
states of $B$ are the new states of the form $p_i$.

If a word $w$ is accepted by the Cho automaton $B$ then $B$ has a
path from the initial state to a some state $p$ labelled $(g, w)$ for
some $g \in S_p$. Since $A$ contains all the states and edges of $B$, we
deduce that $A$ also has a path from the initial state
to $p$ labelled $(g,w)$. Now $g$ must lie in one of the sets $S_p^i$.
It follows easily that $A$ has a path from $p$ to $p_i$ labelled
$(-g, \epsilon)$. Thus, $A$ has a path from the initial state to the
terminal state $p_i$ with label $(1, w)$, so that $A$ accepted $w$.

Conversely, if $A$ accepts a word $w$ then it must have a path from
the initial state to some state $p$ labelled $(g,w)$ and then a path
from $p$ to some terminal state $p_i$ labelled $(-g, \epsilon)$
Clearly the former path must exist also in $B$ and we must have
$g \in S_p^i \subseteq S_p$, so that $w$ is also accepted by $B$.

Thus, the $\Z^k$-automaton $B$ accepts the same language as the
Cho automaton $A$, that is, the word problem for $H$, and so (iii)
holds.

Finally, suppose (iii) holds. Then by
Theorem~\ref{thm_maingauto} and Corollary~\ref{cor_nondetimpdet},
$H$ has word problem accepted by a deterministic $\Z^k$-automaton $A$
for some $k$. We shall construct from $A$ a deterministic Cho automaton
with $2k$ counters recognising the same language.

Suppose $\Z^k$ is generated by
$$\lbrace v_1, \dots, v_k \rbrace$$
and $\N^{2k}$ is generated by
$$\lbrace v_1, \overline{v}_1, \dots, v_k, \overline{v}_k \rbrace.$$
We define a surjective morphism $\sigma : \N^{2k} \to Z^k$ by $v_i \sigma = v_i$ and
$\overline{v}_i \sigma = v_i^{-1}$ for $1 \leq i \leq k$. Define a set
$$S = \lbrace \lambda_1 (v_1 + \overline{v}_1) + \lambda_2 (v_2 + \overline{v}_2) + \dots + \lambda_k (v_k + \overline{v}_k) \mid \lambda_1, \lambda_2, \dots, \lambda_k \in \N \rbrace.$$
Then $S$ is semilinear, and it is easily seen that $S = 1 \sigma^{-1}$.

We now define an $2k$-counter Cho automaton $B$ with the same state set as
$A$, the same designated initial state, and edges constructed as follows.
For each edge in $A$ from $p$ to $q$ with label $(g,w)$, we choose some element $v \in \N^{2k}$ with
$v \sigma = g$, and give $B$ an edge from $p$ to $q$ with label $(v,w)$.
We define $S_q = S$ if $q$ is a terminal state in $A$, and $S_q = \emptyset$
otherwise.

An easy inductive argument shows that $A$ has a path from $p$ to $q$ with
label $(g,w)$ if and only if $B$ has a path from $p$ to $q$ with label
$(v,w)$ for some $v$ such that $v \sigma = g$. In particular, $A$ has a
path from an initial state to a terminal state with label $(1,w)$ if and
only if $B$ has a path from the initial state to a state $q$ with $S_q = S$
having label $(v,w)$ for some $v$ with $v \sigma = 1$, that is, for some
$v \in S$. Hence $w$ is accepted by $A$ if and only if $w$ is accepted by
$B$, as required to complete the proof.
\end{proof}

Note that Theorem~\ref{thm_cho} does not quite provide a complete
characterisation of groups with word problems recognised by $k$-counter
Cho automata for each value of $k$. The proof of Theorem~\ref{thm_cho}
shows only that a group which is virtually free abelian of rank $k$ has word
problem accepted by a $2k$-counter Cho automaton, while a group with word
problem accepted by a $k$-counter Cho automaton is virtually free abelian
of rank $k$ or less. We conjecture that the former bound is tight while
the latter can be strengthened.
\begin{conjecture}
Suppose $H$ is a group with word problem accepted by a $2k$-counter or
$(2k+1)$-counter Cho automaton. Then $H$ is virtually free abelian of
rank $k$.
\end{conjecture}

\section*{Acknowledgements}

This research was started while the second author was at Universit\"at Kassel
supported by a Marie Curie Fellowship within the 6th European Community
Framework Programme; it was completed while he was supported by an RCUK
Academic Fellowship in Manchester. The third author gratefully acknowledges
partial support from European Community funds during a visit to Kassel. The
authors would like to thank Hong-Ray Cho, Sean Cleary, John Fountain, Bob
Gilman, Elaine Render and Benjamin Steinberg for many helpful conversations.

\bibliographystyle{plain}

\def\cprime{$'$} \def\cprime{$'$}

\end{document}